\documentclass[12pt]{article}
\usepackage{amssymb,xcolor,epic}
\usepackage{amsmath}
\usepackage{graphicx}
\usepackage{enumerate}
\usepackage{mathrsfs}
\usepackage{latexsym}
\usepackage{psfrag}
\usepackage{mathrsfs}
\usepackage{multicol}
\usepackage{algorithm}
\usepackage{algpseudocode}

\newcommand{\qed}{\hfill $\square$}

\newenvironment{proof1}{\noindent{\em Proof of Theorem 2.5}.}{\qed\bigskip}
\newenvironment{proof2}{\noindent{\em Proof of Theorem 2.6}.}{\qed\bigskip}

\newtheorem{theorem}{Theorem}[section]
\newtheorem{lemma}[theorem]{Lemma}
\newtheorem{definition}[theorem]{Definition}
\newtheorem{proposition}[theorem]{Proposition}

\newtheorem{conjecture}{Conjecture}[section]

\newcommand{\comments}[1]{}
\begin{document}

\title{Some Sufficient Conditions for Finding a Nesting of the Normalized Matching Posets of Rank 3}

\author{
Yu-Lun Chang
\thanks{Department of Applied Mathematics, National Chung Hsing University, Taichung 40227, Taiwan
{\tt Email:ssn700316@hotmail.com} supported by MOST-105-2115-M-005-003-MY2.}
\and
Wei-Tian Li
\thanks{Department of Applied Mathematics, National Chung Hsing University, Taichung 40227, Taiwan
{\tt Email:weitianli@nchu.edu.tw} supported by MOST-105-2115-M-005-003-MY2.}
}

\date{\small \today}

\maketitle

\begin{abstract}

Given a graded poset $P$, consider a chain decomposition $\mathcal{C}$ of $P$. If $|C_1|\le |C_2|$ implies that the set of the ranks of elements in $C_1$  is 
a subset of the ranks of elements in $C_2$ for any chains $C_1,C_2\in \mathcal{C}$, 
then we say $\mathcal{C}$ is a nested chain decomposition (or nesting, for short) of $P$,  
and $P$ is said to be nested.
In 1970s, Griggs conjectured that every normalized matching rank-unimodal poset is nested.
This conjecture is proved to be true only for all posets of rank 2~\cite{W:05}, some posets of rank 3~\cite{HLS:09,ENSST:11}, 
and the very special cases for higher ranks.
For general cases, it is still widely open. 
In this paper, we  provide some sufficient conditions on the rank numbers of posets of rank 3 to satisfies the Griggs's conjecuture.
\end{abstract}


\section{Introduction}

We start with the necessary terminology of poset theory.
A {\em poset} $P=(P,\le)$ is a set $P$ equipped with a partially order relation $\le$. 
Through out the paper, all posets are finite.
Let $P$ be a poset. We say a subposet $C$ of $P$ is a {\em chain} of length $\ell$ 
if $C=\{x_i\mid x_i< x_{i+1}\mbox{ for }0\le i \le \ell-1\}$.
A {\em chain decomposition} $\mathcal{C}$ of $P$ is a collection of disjoint chains of $P$ with $\cup_{C\in \mathcal{C}}C=P$.
We are looking for decompositions with as few number of chains as possible. 
The most significant theorem in the literature was given by Dilworth~\cite{RPD:1950}. Here an {\em antichain} is a subposet 
of $P$ such that neither $x\le y $ nor $y\le x$ holds for any $x\neq y$ in $Q$.

\begin{theorem}{\rm \cite{RPD:1950}}
For a poset $P$, the minimum number of chains in a chain decomposition is equal to the maximum size of an antichain of $P$.
\end{theorem}

In the following, we study the chain decompositions of a special class of posets, which includes example such as Boolean lattices, linear lattices, and divisor lattices, etc. A {\em graded poset} is a poset such that every maximal chain has the same length. 
For a graded poset $P$, we define the rank function $r:P\longrightarrow \mathbb{N}$ such that 
$r(x)=i$, if there exactly $i$ elements $y< x$ in a maximal chain. 
Moreover, an element $x$ is of rank $i$ if $r(x)=i$ and the rank of $P$ is $\max_{x\in P}r(x)$. 
The {\em $i$th level} of a graded poset $P$ is the collection of all elements of rank $i$, that is, $L_i=\{x\mid r(x)=i, x\in P\}$.
By the definitions, every level is an antichain. 
Therefore, $\max_{i}|L_i|\le |\mathcal{C}|$ for every chain decomposition $\mathcal{C}$ of $P$.
For graded posets, we define a special chain decomposition:
\begin{definition}[nested chain decomposition]
Let $\mathcal{C}$ be a chain decomposition of a graded poset $P$. 
For any chains $C_i,C_j\in \mathcal{C}$, if $|C_i|\le |C_j|$ implies 
$\{r(x)\mid x\in C_i\}\subseteq \{r(x)\mid x\in C_j\}$, then $\mathcal{C}$ is called a 
nested chain decomposition of $P$. 
We say $P$ is nested, or it has a nesting, if such a decomposition of $P$ exists.
\end{definition}
Observe that a nesting $\mathcal{C}$  is a chain decomposition with minimum number of chains. 
Because from the inclusion relation, $\cap_{C\in \mathcal{C}}\{r(x)\mid x\in C\}$ is not empty, 
and there exists some $m$ and a level $L_m$
such that $m\in\cap_{C\in\mathcal{C}}\{r(x)\mid x\in C\}$. 
Thus, every chain in $\mathcal{C}$ contains an element of rank $m$, and hence $|L_m|=|\mathcal{C}|$. 
Since $\max_i|L_i|\le |\mathcal{C}|=|L_m|\le \max_i|L_i|$, we have $\max_i|L_i|\le |\mathcal{C}|$.
We refer the reader to see more properties of graded posets in~\cite{A:02,ENG:97}. 

Anderson\cite{A:1967b} and Griggs\cite{JRG:1977} independently gave the same sufficient condition for the existence of a nesting in the graded posets. 
Let $r_i$ denote the cardinality of level $i$. The {\em rank numbers} or {\em Whitney numbers} of a graded poset of rank $n$ is the 
sequence $(r_0,r_1,\ldots, r_n)$. If $r_i=r_{n-i}$ for all $0\le i\le n$, then $P$ is said to be {\em rank-symmetric}. 
Suppose there exists some $i$ such that $r_0\le\cdots \le r_i\ge \cdots \ge r_n$, then $P$ is {\em rank-unimodal}.
For any levels $L_i$ and $L_j$ of $P$, consider any subset $S$ of $L_i$ and denote the set 
$\Gamma_j(S)=\{x\in L_j\mid x \le y\mbox { or }y\le x \mbox{ for some  }y\in S\}$. 
If the inequality
\[\frac{|S|}{r_i} \le \frac{|\Gamma_j(S)|}{r_j}\]
holds, then we say $P$ has {\em the normalized matching property} from $i$ to $j$. 
By simple calculation, one can see if $P$ has the normalized matching property from $i$ to $j$, 
then it also has the property from $j$ to $i$.
Moreover, if $P$ has the normalized matching property from $i$ to $j$ and $j$ to $k$ for $i<j<k$, then 
it has the property from $i$ to $k$.
Once the normalized matching property holds between any two levels, then we say $P$ is a {\em normalized matching poset}. 

\begin{theorem}{\rm \cite{A:1967b,JRG:1977}}\label{thm:scd}
A graded poset $P$ has a nesting if it is a normalized matching poset, and is rank-symmetric and rank-unimodal.
\end{theorem}

In fact, every chain in a nesting of a poset $P$ described in the theorem 
contains elements of ranks $i,i+1,\ldots, n-i$ for some $i$, where $n$ is the rank of $P$.   
Such a decomposition is also called a {\em symmetric chain decomposition} of $P$. 
In addition to Theorem~\ref{thm:scd}, Griggs also posed the following conjecture~\cite{JRG:1977,JRG:1988,JRG:1995}:

\begin{conjecture}[\rm Griggs Nesting Conjecture]\label{gnc:un}
Every normalized matching rank-unimodal poset is nested.
\end{conjecture}

The conjecture turns out to be extremely difficult, although one can easily give an affirmative answer of graded posets of rank 1 using 
the well-known Hall's Marriage Theorem~\cite{PH:1935}. 
There are only a few graded posets of small ranks which are proven to satisfy the conjecture by Wang~\cite{W:05}, 
Hsu, Logan, and Shahriari.~\cite{HLS:09}, and 
Escamilla, Nicolae, Salerno, Shahriari, and Tirrell~\cite{ENSST:11}, respectively. 
In section 2, we introduce the early results on the graded posets of ranks 2 and 3, 
and mention our theorems at the end of the section. 
The main contribution of this paper is to give more sufficient conditions on the graded posets of rank 3 to satisfy the conjecture,
based on the ideas of proofs in the early papers.  
The proofs of our theorems are presented in Section 3.

\section{Normalized Matching Posets of Rank 2 and 3}

Note that in Conjecture~$\ref{gnc:un}$, we only concern the normalized matching property and the conditions of the rank numbers. 
The structure of the poset is irrelevant. 
For convenience, we use the notation $NM(r_0,r_1,\ldots,r_n)$ to denote the collection 
of all normalized matching posets of rank $n$ with rank numbers $(r_0,r_1,\ldots,r_n)$. 
In 2005, Wang~\cite{W:05} dropped the rank-unimodal assumption and proved a stronger result. 

\begin{theorem}{\rm \cite{W:05}}\label{thm:05}
Every poset $P\in NM(r_0,r_1,r_2)$ has a nesting. 
\end{theorem}

For graded posets of rank 3, Shahriari with two research groups\cite{HLS:09,ENSST:11} developed some sufficient conditions on 
the rank numbers $(r_0,r_1,r_2,r_3)$ to guarantee the existence of a nesting. 
In~\cite{HLS:09}, the authors came up with a clever idea which can not only simplify the proof of Theorem~\ref{thm:05} but also reduce 
the rank numbers to fewer cases that need to be considered for graded posets of rank 3.
Since we will use this idea in our proof, we introduce it below.
\begin{proposition}
Given a graded poset $P$, let $P'=P\cup \{x\}$ be obtained by adding a new element $x$ to the $i$th level of $P$ together
with the partial order relations $y< x$ (resp. $x<y$) if $y\in L_j$ and $j<i$ (resp. $j<i$).  
If $P\in NM(r_0,r_1,\ldots, r_n)$, then $P'\in NM(r_0,\ldots,r_{i-1},r_{i}+1,r_{i+1},\ldots, r_n)$. 
\end{proposition} 
The proof of the proposition is straightforward, since if we pick a set $S$ in the $i$th level of $P'$, 
either it contains $x$, then $|\Gamma_k(S)|/r_k=1$, or it does not contain $S$, then $|S|/(r_i+1)< |S|/r_i< |\Gamma_k(S)|/r_k$, 
for all $k$. In~\cite{HLS:09}, such an element is called a {\em ghost}.

We demonstrate two instances of exploiting the ghosts to get a nesting. 
Suppose $P$ is a poset in $NM(r_0,r_1,r_2)$ with $r_0<r_2<r_1$. 
Then we add $r_2-r_0$ ghosts to the $0$th level to get a rank-symmetric poset $P'$. 
By Theorem~\ref{thm:scd}, $P'$ has a nesting and each chain contains elements of ranks either $\{0,1,2\}$ or $\{2\}$. 
After removing the ghosts from the chains of length 2, 
we obtain a chain decomposition of $P$ such that each chain contains elements of ranks either $\{0,1,2\}$, or $\{1,2\}$, or $\{2\}$.
If the rank numbers satisfy $r_1<r_0<r_2$, then we add $r_0-r_1$ ghosts to the first level. 
The new poset $P''$ restricted on the first and second levels is a poset of rank 1. 
So we can partition it into chains of length either 0 or 1. 
Meanwhile, the poset consisting of the 0th and first level of $P''$ can be partitioned 
into chains of length of 1. 
The chains of length 1 in two decompositions can be concatenated into chains of length 2.
Finally, we remove the ghosts to get of decomposition of $P$ with each chain containing elements of 
rank either $\{0,1,2\}$, or $\{0,2\}$, or $\{2\}$.
Indeed, the arguments above are exact the ideas of Hsu et al. in~\cite{HLS:09}, 
used to reprove Theorem~\ref{thm:05} .

Using the ideas of the ghost elements, the induction, and the duality, 
Hsu et al.~\cite{HLS:09} showed that to prove Conjecture~$\ref{gnc:un}$ for posets of rank 3, it suffices to verify that all posets $P\in (r_0,r_1,r_2,r_3)$ with $r_2>r_1>r_0=r_3$ are nested.
For example, if a poset $P\in (r_0,r_1,r_2,r_3)$ with $r_2>r_1>r_0>r_3$, 
then we add $r_0-r_3$ ghosts to the third level of $P$ to get a new poset $P'\in (r_0,r_1,r_2,r_0)$. 
Now suppose we already have a nesting of $P'$. 
We then remove the ghosts in all the longest chains to get a nesting of $P$.
See~\cite{HLS:09} for the details of all the reduction methods.
With this assumption on the rank numbers, Hsu et al.~\cite{HLS:09} 
and Escamilla et al.~\cite{ENSST:11} proved the following Theorem~\ref{thm:09} and Theorem~\ref{thm:11}, respectively.

\begin{theorem}{\rm \cite{HLS:09}}\label{thm:09}
Let $r_0, r_1$ and $r_2$ be positive integers with $r_0 <r_1 <r_2$. Assume that at least one of the following conditions are satisfied:
\begin{description}
\item{(a)}   $r_{1} \ge r_{2}-\lceil{\frac{r_{2}}{r_{0}}}\rceil+1$;
\item{(b)}   $r_{1}$ $=r_{0}+1$;
\item{(c)}   $r_{2}> r_{0}r_{1}$;
\item{(d)}   $r_{1}$ divides $r_{2}$.
\end{description}
Then every $P\in NM(r_0,r_1,r_2,r_0)$ is nested. 
\end{theorem}

\begin{theorem}{\rm \cite{ENSST:11}}\label{thm:11}
Let $r_0, r_1$ and $r_2$ be positive integers with $r_0 <r_1 <r_2$. Assume that at least one of the following conditions are satisfied:
\begin{description}
\item{(a)} $r_{0}$ divides $r_{1}$, or
\item{(b)} $r_{0}+1$ divides $r_{1}$, or
\item{(c)}  $f(i)\geq 0$ for all $1\leq i\leq r_{1}-r_{0}-1$, where the function $f$ is defined by 
\[
f(i) = \left\lceil\frac{r_{0}(1+i)}{r_{2}-r_{0}}\right\rceil - 
\left\lfloor\frac{r_{0}i}{r_{1}-r_{0}}\right\rfloor\mbox{; or}\] 
\item{(d)} $r_{2}> r_{0}r_{1}-r_{0}\gcd(r_{1},r_{2})$.
\end{description}
Then every $P\in NM(r_0,r_1,r_2,r_0)$ is nested. 
\end{theorem} 

In~\cite{ENSST:11}, the authors also examined the posets of rank 3 with $r_2\le 13$. 
Using Theorem~\ref{thm:11}, one can verify that if $r_0<r_1<r_2\le 13$, every poset
$P\in NM(r_0,r_1,r_2,r_0)$ satisfies Conjecture~\ref{gnc:un} except that 
$(r_0,r_1,r_2,r_0)$ is equal to one of the six cases:
$(6,8,12,6)$, $(6,9,12,6)$, $(4,6,13,4)$, $(5,8,13,5)$, $(6,8,13,6)$, $(6,9,13,6)$. 
We close Section 2 by stating our results.
For graded posets of rank 3, 
we provide two more sufficient conditions on the rank numbers for the existence of a nesting:

\begin{theorem}\label{thm1} 
Let $P\in NM(r_0,r_1,r_2,r_0)$.
If both $r_{1}$ and $r_{0}$ divide $r_{2}-1$, then $P$ is nested. 
\end{theorem}

\begin{theorem}\label{thm2}
Let $P\in NM(r_0,r_1,r_2,r_0)$.
If $kr_{0}\le r_{1}\le k(r_{0}+1)$, then $P$ is nested.
\end{theorem}

Observe that by Theorem~\ref{thm1}, we see that every $P\in NM(4,6,13,4)$ has a nesting. 
Unfortunately, other unsolved cases with $r_2\le 13$ mentioned in\cite{HLS:09} cannot be settled by our theorems. 
Nevertheless, for $14\le r_2\le 15$, we can use Theorem~\ref{thm2} to show every $P\in NM(r_0,r_1,r_2,r_3)$ with
$(r_0,r_1,r_2,r_3)\in\{(4,9,14,4)$, $(3,7,15,3)$,$(4,9,15,4)$,$(5,11,15,5)\}$ is nested. 
These are not covered by Theorem~\ref{thm:09} and~\ref{thm:11}.

\section{Proofs of the Main Theorems}
In this section, we give the proofs of Theorem~\ref{thm1} and Theorem~\ref{thm2} .
Let us begin with the proof of Theorem~\ref{thm1}.

\begin{proof1}
Pick a poset $P\in NM(r_0,r_1,r_2,r_0)$, where
$k_0r_0+1=r_2$ and $k_1r_1+1=r_2$ for some integers $k_1$ and $k_2$. 
Let $P'$ be a poset obtained by removing an arbitrary element $x$ from $L_2$ .
To show that $P'$ is a normalized matching poset, 
we only need to verify the inequality holds between $L_1$ and $L_2\setminus\{x\}$ as well as $L_2\setminus\{x\}$ and $L_3$. 
First consider $L_1$ and $L_2\setminus\{x\}$.
By the symmetry, we only need to verify the normalized matching property from $L_1$ to $L_2\setminus \{x\}$.
For any $S\subseteq L_1$, since $P$ is a normalized matching poset, we have
\[
\frac{|S|}{r_1}\le \frac{|\Gamma_2(S)|}{r_2}=\frac{|\Gamma_2(S)|}{k_1r_1+1}.
\]
Equivalently,
\[{k_{1}|S|}+\frac{|S|}{r_1}\le |\Gamma_2(S)|.\] 
When $S\neq\emptyset$, we have $k_1|S|+1\le |\Gamma_2(S)|$ since $|\Gamma_2(S)|$ is an integer. 

Now, for $P'$, if the removed element  $x$ is not in $\Gamma_2(S)$, 
then 
\[
\frac{|S|}{r_{1}}\le \frac{|\Gamma_2(S)|}{k_1r_1+1}\le \frac{|\Gamma_2(S)|}{k_1r_1}=\frac{|\Gamma_2(S)|}{r_{2}-1}=\frac{|\Gamma_2(S)|}{|L_2\setminus \{x\}|}.
\] 
Otherwise,  $x\in \Gamma_2(S)$ and then $S\neq\emptyset$. 
We have 
\[
\frac{|S|}{r_{1}}=\frac{k_{1}|S|}{k_{1}r_{1}}\le \frac{|\Gamma_2(S)|-1}{r_2-1}=\frac{|\Gamma_2(S)|-1}{|L_2\setminus \{x\}|}.
\]
The numerator $|\Gamma_2(S)|-1$ is just the number of elements $y\in L_2\setminus\{x\}$ satisfying $z<y$ for some $z\in S$.
So the normalized matching property holds between $L_1$ and $L_2\setminus\{x\}$. 
Using a similar argument we can see that the normalized matching property also holds between $L_3$ and $L_2\setminus\{x\}$. 
Now that $r_1$ divides $k_1r_1=r_2-1$, 
so $P'$ has a nesting $\mathcal{C}$ by Theorem~\ref{thm:09} (d). 
Finally, we view $\{x\}$ as a one-element chain and add it to $\mathcal{C}$ to get a nesting of $P$.         
\end{proof1}

It is worth mentioning that the proof in Theorem~\ref{thm1} is similar to the next lemma in~\cite{HLS:09},  
which is used to prove Theorem~\ref{thm:09} (b). 

\begin{lemma}{\rm \cite{HLS:09}}\label{avoid}
Let $P\in NM(r_0,r_0+1)$. For any $x$ of rank 1 in $P$, there exists a chain partition of $P$ which consists of  
$r_0$ chains of length 1 and another chain $\{x\}$ of length 0. 
\end{lemma}

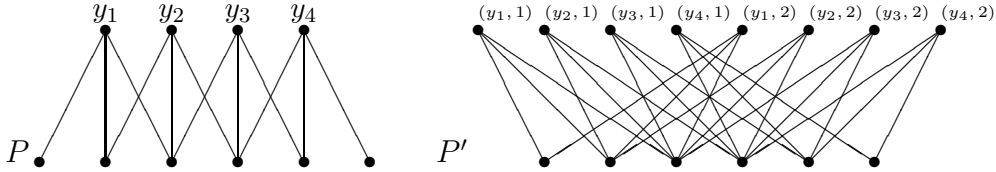
\begin{figure}[h]
\begin{center}
$P$
\begin{picture}(150,60)
\put(0,0){\circle*{4}}
\put(25,0){\circle*{4}}
\put(50,0){\circle*{4}}
\put(75,0){\circle*{4}}
\put(100,0){\circle*{4}}
\put(125,0){\circle*{4}}

\put(25,50){\circle*{4}}
\put(50,50){\circle*{4}}
\put(75,50){\circle*{4}}
\put(100,50){\circle*{4}}

\put(20,55){{\small $y_1$}}
\put(45,55){{\small $y_2$}}
\put(70,55){{\small $y_3$}}
\put(95,55){{\small$y_4$}}

\put(25,50){\line(-1,-2){25}}
\put(50,50){\line(-1,-2){25}}
\put(75,50){\line(-1,-2){25}}
\put(100,50){\line(-1,-2){25}}

\put(25,50){\line(0,-1){50}}
\put(50,50){\line(0,-1){50}}
\put(75,50){\line(0,-1){50}}
\put(100,50){\line(0,-1){50}}

\put(25,50){\line(1,-2){25}}
\put(50,50){\line(1,-2){25}}
\put(75,50){\line(1,-2){25}}
\put(100,50){\line(1,-2){25}}
\end{picture}$P'$
\begin{picture}(200,60)
\put(0,50){\circle*{4}}
\put(25,50){\circle*{4}}
\put(50,50){\circle*{4}}
\put(75,50){\circle*{4}}
\put(100,50){\circle*{4}}
\put(125,50){\circle*{4}}
\put(150,50){\circle*{4}}
\put(175,50){\circle*{4}}

\put(25,0){\circle*{4}}
\put(50,0){\circle*{4}}
\put(75,0){\circle*{4}}
\put(100,0){\circle*{4}}
\put(125,0){\circle*{4}}
\put(150,0){\circle*{4}}

\put(0,50){\line(1,-2){25}}
\put(0,50){\line(1,-1){50}}
\put(0,50){\line(3,-2){75}}
\put(25,50){\line(1,-2){25}}
\put(25,50){\line(1,-1){50}}
\put(25,50){\line(3,-2){75}}
\put(50,50){\line(1,-2){25}}
\put(50,50){\line(1,-1){50}}
\put(50,50){\line(3,-2){75}}
\put(75,50){\line(1,-2){25}}
\put(75,50){\line(1,-1){50}}
\put(75,50){\line(3,-2){75}}

\put(100,50){\line(-1,-2){25}}
\put(100,50){\line(-1,-1){50}}
\put(100,50){\line(-3,-2){75}}
\put(125,50){\line(-1,-2){25}}
\put(125,50){\line(-1,-1){50}}
\put(125,50){\line(-3,-2){75}}
\put(150,50){\line(-1,-2){25}}
\put(150,50){\line(-1,-1){50}}
\put(150,50){\line(-3,-2){75}}
\put(175,50){\line(-1,-2){25}}
\put(175,50){\line(-1,-1){50}}
\put(175,50){\line(-3,-2){75}}

\put(0,55){{\tiny $(y_1,1)$}}
\put(25,55){{\tiny $(y_2,1)$}}
\put(50,55){{\tiny $(y_3,1)$}}
\put(75,55){{\tiny $(y_4,1)$}}

\put(100,55){{\tiny $(y_1,2)$}}
\put(125,55){{\tiny $(y_2,2)$}}
\put(150,55){{\tiny $(y_3,2)$}}
\put(175,55){{\tiny $(y_4,2)$}}
\end{picture}
\end{center}
\caption{A new poset obtained by replacing a 2-clone of the first level of $P$ to it.}\label{fig1}
\end{figure}

Before presenting the proof of Theorem~\ref{thm2}, we need more preparations.
In addition to adding the ghosts to a normalized matching poset, 
there are some techniques to produce new normalized matching posets from the old ones. 
We introduce two construction approaches.
\begin{definition}[$k$-clone]
Let $P$ be a graded poset and $L_i$ be a level of $P$. Then $L$ is said to be a $k$-clone of $L_i$ if $L=L_i\times\{1,\ldots,k\}$ and the partial 
order relations of each $(y,i)$ and others elements in $L_{i-1}$ (resp. $L_{i+1}$) is $x<(y,i)$ (resp. $(y,i)<z$) if and only 
if there exist some $x\in L_{i-1}$ and $y\in L_i$ (resp. $y\in L_j$ and $z\in L_{i+1}$). See Figure~\ref{fig1} as an illustration.
\end{definition}

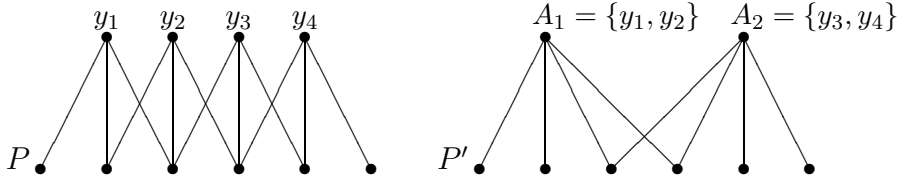
\begin{figure}[h]
\begin{center}
$P$
\begin{picture}(150,60)
\put(0,0){\circle*{4}}
\put(25,0){\circle*{4}}
\put(50,0){\circle*{4}}
\put(75,0){\circle*{4}}
\put(100,0){\circle*{4}}
\put(125,0){\circle*{4}}

\put(25,50){\circle*{4}}
\put(50,50){\circle*{4}}
\put(75,50){\circle*{4}}
\put(100,50){\circle*{4}}

\put(20,55){{\small $y_1$}}
\put(45,55){{\small $y_2$}}
\put(70,55){{\small $y_3$}}
\put(95,55){{\small$y_4$}}

\put(25,50){\line(-1,-2){25}}
\put(50,50){\line(-1,-2){25}}
\put(75,50){\line(-1,-2){25}}
\put(100,50){\line(-1,-2){25}}

\put(25,50){\line(0,-1){50}}
\put(50,50){\line(0,-1){50}}
\put(75,50){\line(0,-1){50}}
\put(100,50){\line(0,-1){50}}

\put(25,50){\line(1,-2){25}}
\put(50,50){\line(1,-2){25}}
\put(75,50){\line(1,-2){25}}
\put(100,50){\line(1,-2){25}}
\end{picture}$P'$
\begin{picture}(150,60)
\put(0,0){\circle*{4}}
\put(25,0){\circle*{4}}
\put(50,0){\circle*{4}}
\put(75,0){\circle*{4}}
\put(100,0){\circle*{4}}
\put(125,0){\circle*{4}}

\put(25,50){\circle*{4}}
\put(100,50){\circle*{4}}

\put(20,55){{\small $A_1=\{y_1,y_2\}$}}
\put(95,55){{\small$A_2=\{y_3,y_4\}$}}

\put(25,50){\line(-1,-2){25}}
\put(100,50){\line(-1,-2){25}}
\put(25,50){\line(0,-1){50}}
\put(100,50){\line(0,-1){50}}
\put(25,50){\line(1,-2){25}}
\put(100,50){\line(1,-2){25}}
\put(25,50){\line(1,-1){50}}
\put(100,50){\line(-1,-1){50}}
\end{picture}
\end{center}
\caption{A new poset obtained by replacing a 2-bunch of the first level of $P$ to it.}\label{fig2}
\end{figure}

\begin{definition}[$m$-bunch] 
Let $P$ be a graded poset and $L_i$ be a level of $P$. Suppose $|L_i|=m\ell$ for some integers $m$ and $\ell$.
First partition $L_1$ into $m$ arbitrary subsets $A_{1}$,\ldots,$A_{m}$ of equal size $\ell$.
Then $L$ is an $m$-bunch of $L_i$ if $L=\{A_{1}$,\ldots,$A_{m}\}$ and the partial 
order relations of each $A_j$ and others elements in $L_{i-1}$ (resp. $L_{i+1}$) is $x<A_j$ (resp. $A_j<z$)
if and only if there exist $x\in L_{i-1}$ and $y\in A_j$ (resp. $y\in A_j$ and $z\in L_{i+1}$). 
See Figure~\ref{fig2} as an illustration.
\end{definition}

The above operations on posets preserve the normalized matching property:
\begin{proposition}\label{prop:CB}
If $P$ is a normalized matching poset, 
then the new poset obtained by replacing a $k$-clone or an $m$-bunch of some level
of $P$ to it is still a normalized matching poset.
\end{proposition}
The proof of this proposition was given by Hsu et al.~\cite{HLS:09} (clone),
and by Escamilla et al.~\cite{ENSST:11} (bunch), respectively.
Now we prove our second theorem.

\begin{proof2}
Consider $P\in NM(r_0,r_1,r_2,r_0)$ with $kr_0\le r_1\le k(r_0+1)$ for some integer $k$. 
Note that the two ends of the inequality are in the statements of Theorem~\ref{thm:11} (a) and (b). 
Thus, we may suppose $r_{1}=kr_{0}+t$ for some $1\le t\le k-1$.
Pick a poset $P\in NM(r_0,r_1,r_2,r_0)$. 
We use the induction method to find the nestings of subposets induced by different levels of $P$. 
Our goal is to combine the nestings properly to get a nesting of $P$.

First construct a poset $P_1$ of rank 2 induced by the top three levels of $P$ 
with a replacement of a $k$-clone of the highest level.  
By Proposition~\ref{prop:CB}, $P_1\in NM(r_{1},r_{2},kr_{0})$, 
and there exists a nesting $\mathcal{C}_1$ of $P_1$ by Theorem~\ref{thm:05}. 
Observe that there are $kr_0$ chains of length 2 in $\mathcal{C}_1$ such that eahc of them contains an 
element $(y,j)$ in the highest level of $P_1$.

Clearly, the bottom two levels $L_0$ and $L_1$ of $P$ induce a subposet of rank 1 and has a nesting.  
However, we do not want a nesting of the above poset containing a chain of length 1 whose top element
is the bottom element of a chain of length 1 in $\mathcal{C}_1$. 
This could lead to two chains of length 2 but the ranks of elements in one chain is $\{0,1,2\}$ and the other is $\{1,2,3\}$
when we combine the two nestings together.
To avoid this, we construct a poset $P_2$ of rank 1 as follows.
At the beginning, we add $k-t$ additional ghosts to $L_1$ of $P$ in advance.
Now this level contains $k(r_0+1)$ elements, and we will partition them into $r_0+1$ sets of size $k$.
Because $L_1$ is also the bottom level of $P_1$, for each $y_i\in L_3$ of $P$, there exist exactly $k$ elements in $L_1$ such that each of them lies in a chain, 
containing $(y_i,j)$ for some $1\le j\le k$, of length 2 in $\mathcal{C}$.
In addition, there are $t=r_1-kr_0$ elements in $L_1$ which are not in any chain of length 2 in $\mathcal{C}$. 
For $1\le i\le r_{0}$, let $A_{i}$ be the set consisting of every element in $L_1$, which lies in a chain in $\mathcal{C}_1$ containing the element $(y_i, j)$ for some $1\le j\le k$ .
Moreover, let $A_{r_{0}+1}$ be the set consisting of the $t$ remaining elements in $L_1$ and the $k-t$ ghosts. 
We bunch all elements in $L_1$ and the ghosts into the above $r_0+1$ sets $A_1,\ldots, A_{r_0+1}$.
The poset induced by these $A_i$s and $L_0$ is $P_2$.

By Lemma $\ref{avoid}$, there is a chain partition $\mathcal{C}_2$ of $P_2$ 
with $r_0$ chains of length 1 and one chain of length 0
such that each chain of length 1 does not contain $A_{r_0+1}$ and each $x_i\in L_0$ is in a chain of length 1 in $\mathcal{C}_2$.
Assume the chains are $\{x_i,A_i\}$ for $1\le i\le r_0$. 
It follows that for each $i$ there exists some element in $z\in A_i$ with $x_i<z$. 
Fix some $i$. For those $k$ chains of length 2 in $\mathcal{C}_1$ containing $(y_i,j)$ for some $1\le j\le k$, 
we extend one of them to length 3 by adding the element $x_i$
and delete the top elements of the remaining $k-1$ chains of length 2.
Repeating the operations for all $1\le i\le r_0$ gives us a nesting of $P$.
\end{proof2}

\end{document}